\documentclass[10pt,twoside]{article}
\usepackage{amssymb}
\usepackage{fullpage}
\usepackage{a4,amsmath}
\usepackage{amsmath}
\usepackage{amscd}

\def\R{\mathbb{R}}

\def\1{\mathbf{1}}
\def\:{\lrcorner}
\def\#{\sharp}

\def\e{\epsilon}

\def\x{\otimes}
\def\<#1,#2>{\langle#1,\,#2\rangle}

\def\qed{\ensuremath{\quad\Box\quad}}
\def\pfill{\par\vskip2mm plus1mm minus1mm\noindent}
\def\inv#1{\raise.1em\hbox to 0pt{$^{-1}$\hss}_{#1}\;}

\def\v{\noindent}

\newcommand{\bean}{\begin{eqnarray*}}
\newcommand{\eean}{\end{eqnarray*}}
\newcommand{\benu}{\begin{enumerate}}
\newcommand{\eenu}{\end{enumerate}}
\newcommand{\eea}{\end{eqnarray}}
\newcommand{\bea}{\begin{eqnarray}}
\newtheorem{Theorem}{Theorem}[section]

\title{A note on closed isometric embeddings}
\begin{document}

\author{Olaf M\"uller\footnote{Instituto de Matem\'aticas, Universidad Nacional Aut\'onoma de M\'exico (UNAM) Campus Morelia, C. P. 58190, Morelia, Michoac\'an, Mexico. email: olaf@matmor.unam.mx}}

\maketitle

\begin{abstract}
\v A famous theorem due to Nash (\cite{jN}) assures that every Riemannian manifold can be embedded isometrically into some Euclidean space $E^n$. An interesting question is whether for a {\em complete} manifold $M$ we can find a {\em closed} isometric embedding. This note gives the affirmative answer to this question asked to the author by Paolo Piccione. 
\end{abstract}

\v In his famous 1956 article John Nash proved that every Riemannian metric on an $n$-dimensional manifold $M$ can be constructed as a pullback metric for an embedding of $M$ into some Euclidean space. He gave also an estimate of the smallest possible dimension $N$ of the Euclidean space as $N = \frac{1}{2} \cdot  n \cdot (n+1) \cdot (3n+11)$. Now one can try to find some stronger derivates of this theorem if strengthening the assumptions. In this note, we want to examine the question whether every {\em complete} manifold admits a {\em closed} isometric embedding. Although folk wisdom apparently has a positive answer to this question already, there does not seem to be any proof in the literature up to now. The question is more difficult than it might seem at first sight as there are plenty of non-closed isometric embeddings of complete manifols, e.g. spirals converging to $0$ or to a circle as isometric embeddings of $\R$. 

\bigskip

\v For a Lipschitz function $f$ on a metric space $M$, we denote by $L(f):= \sup \{ \frac{\vert f(p) - f(q) }{d(p, q)} \vert$

\v $  p, q \in M, p \neq q \}$ its Lipschitz number. We will need the following nice theorem from \cite{AFLR} about approximation of Lipschitz functions by smooth functions (even on infinite-dimensional Riemannian manifolds):

\begin{Theorem}
Let $(M, g)$ be a separable Riemannian manifold, let $f: M \rightarrow \R$ be a Lipschitz function, let $\rho: M \rightarrow (0, \infty)$ be a continuous function, and let $r>0$. Then there is a $C^{\infty}$ and Lipschitz function $g:M \rightarrow \R$ with $\vert f(p) - g(p) \vert \leq \rho (p) $ for every $p \in M$, and $L (g) \leq L (f) + r$. \hfill \qed
\end{Theorem}

\v Now let us state and prove our theorem. The basic idea of the proof is to look at balls of increasing radius and to define an imbedding which lifts the larger and larger balls into an additional direction thereby resolving a possible spiralling. As the distance itself is not differentiable in general, we have to be a little bit more careful and thus we will need the theorem above. 

\begin{Theorem}
If $(M,g)$ is a complete $n$-dimensional Riemannian manifold, then there is a closed isometric $C^{\infty}$-embedding of $(M,g)$ into $E^{N+1}$, where $N:= \frac{1}{2} \cdot n \cdot (n+1) \cdot (3n+11)$.
\end{Theorem}

\v {\bf Proof.} Without restriction of generality, let $M$ be noncompact. Choose a point $p \in M$ and define $D: M \rightarrow \R$ by $D(q) := d(p,q)$. Obviously this is a continuous function on $M$, but in general it is not $C^1$ because of possible cut points. Define $f:= \frac{2}{3} D $ outside a small geodesic ball $B_r (p)$ with $r < 1/4$ and $f:= \frac{2}{3} r$ in $B_r (p)$. By some case distinctions and the triangle inequality it is easy to see that this is a Lipschitz function with Lipschitz number $2/3$. Then let $\rho $ as in the previous theorem be determined by the tube $(\frac{3}{4} f, \frac{3}{2} f)$ around $f$, thus we can we find a smooth function $\phi$ contained in this neighborhood. As $f$ is Lipschitz with Lipschitz number $2/3$, we can choose $\phi $ Lipschitz with Lipschitz number $3/4$. We have $D/2 < \phi < D$ outside $B_r (p)$. Now we define a new metric $\tilde{g} := g - \frac{1}{4} d \phi \x d \phi$. Because of the Lipschitz condition, we have $ \vert \vert d \phi \vert \vert < 3/4$, thus $\tilde{g}$ is a smooth Riemannian metric on $M$. With Nash's embedding theorem we find an isometric embedding $\tilde{\e}$ of $(M, \tilde{g})$ into the Euclidean space $E^N$. If we modify this embedding by adding one dimension and defining $\e := \tilde{\e} + \frac{1}{2} \phi \cdot e_{n+1}$, then $\e$ is an isometric embedding for the original metric $g$. Now let a point $ q \in E^{N+1} \setminus \e (M)$ be given; denote by $Q:= q_{N+1}$ its last coordinate. Then on $K:= M \setminus B_{8 Q} (p)$, $\phi$ is greater than $4Q$, thus $K$ is mapped to $E^N  \times [2 Q, \infty ) \subset E^{N+1}$ by $\e$, therefore $d(q, \e (M \setminus \overline{B}_{8 Q} (p))) \geq Q$. On the other hand, because of completeness of $(M, g)$, $\overline{B}_{8Q} (p)$ is compact. Therefore $\e (\overline{B}_{8Q}(p))$ is compact, too, and has a nonzero distance to $q$ as well, thus $\e$ is closed. \hfill \qed

\bigskip

\v The author wants to thank Paolo Piccione for useful comments on the first version of this note.

\end{document}